\documentclass[a4paper,12pt,reqno]{amsart}

\usepackage[indent=1em,skip=0.5em]{parskip}

\usepackage[utf8]{inputenc}
\usepackage{microtype}
\usepackage{amsmath,amsthm,amsfonts,amssymb}
\usepackage{graphicx}
\usepackage[margin=0.65in]{geometry}

\usepackage{xparse}

\usepackage[dvipsnames]{xcolor}
\usepackage[hyphens]{url}
\usepackage{hyperref}
\hypersetup{
    colorlinks=true,
    linkcolor=cyan!80!black,
    citecolor=MidnightBlue,
    urlcolor=magenta,
}
\usepackage[hyphenbreaks]{breakurl}
\usepackage{quoting}

\usepackage{accents}
\renewcommand{\mathring}[1]{\accentset{\circ}{#1}}

\theoremstyle{plain}

\newtheorem{lemma}{Lemma}[section]
\newtheorem{theorem}{Theorem}[section]

\theoremstyle{definition}
\newtheorem{defn}{Definition}[section]

\theoremstyle{remark}
\newtheorem{remark}{Remark}[section]

\newcommand{\A}{{\mathcal{A}}}

\newcommand{\C}{\mathbb{C}}

\newcommand{\one}{{\bf 1}}
\newcommand{\zero}{{\bf 0}}
\newcommand{\NC}{\mathsf{NC}}

\NewDocumentCommand{\Pure}{O{a} m}{{\mathring{#1}}_{#2}}
\NewDocumentCommand{\PureStar}{O{a} m}{({\mathring{#1}}_{#2})^*}
\NewDocumentCommand{\Alg}{O{\mathcal{A}} m}{{\mathring{#1}}_{#2}}
\NewDocumentCommand{\Mixed}{O{a} O{i} m m}{{\mathring{#1}}_{{#2}_{#3}} \otimes \cdots \otimes {\mathring{#1}}_{{#2}_{#4}}}
\NewDocumentCommand{\MixedTwo}{O{a} O{i}}{{\mathring{#1}}_{{#2}_{1}} \otimes {\mathring{#1}}_{{#2}_{2}}}
\NewDocumentCommand{\MixedStar}{O{a} O{i} m m}{({\mathring{#1}}_{{#2}_{#4}})^{*} \otimes \cdots \otimes ({\mathring{#1}}_{{#2}_{#3}})^{*}}
\NewDocumentCommand{\MixedStarTwo}{O{a} O{i}}{({\mathring{#1}}_{{#2}_{2}})^{*} \otimes ({\mathrm{#1}}_{{#2}_{1}})^{*}}
\NewDocumentCommand{\AlgMixed}{O{\mathcal{A}} O{i} m m}{{\mathring{#1}}_{{#2}_{#3}} \otimes \cdots \otimes {\mathring{#1}}_{{#2}_{#4}}}

\numberwithin{equation}{section}

\title[Construction of product $*$-probability spaces via free cumulants]{Construction of product $*$-probability spaces \\ via free cumulants}

\author[A. Bose]{Arup Bose}
\address{Statistics and Mathematics Unit, Indian Statistical Institute, 203 B. T.~Road, Kolkata 700108, India}
\email{bosearu@gmail.com}
\thanks{AB was partially supported by his J.C. Bose Fellowship JBR/2023/000023 from ANRF, Govt.~of India.}

\author[S.S. Mukherjee]{Soumendu Sundar Mukherjee}
\address{Statistics and Mathematics Unit, Indian Statistical Institute, 203 B. T.~Road, Kolkata 700108, India}
\email{ssmukherjee@isical.ac.in}
\thanks{SSM was partially supported by an INSPIRE research grant (DST/INSPIRE/04/2018/002193) from the Dept.~of Science and Technology, Govt.~of India and a Start-Up Grant from Indian Statistical Institute, Kolkata.}


\begin{document}

\maketitle
\begin{quoting}[vskip=3pt,font={itshape}]
    \scriptsize
    Professor K.R. Parthasarathy (KRP in the Stat-Math ISI tradition) is a revered name across the world of Probability. A master of his craft, he was always full of energy, and ever willing to share his knowledge in his uniquely captivating style with anyone who would be interested. The first author fondly recalls the encouragement he received as a young student and faculty from KRP on multiple occasions, and the long one-on-one conversation he had with him a few years ago. Although the second author never had the opportunity to meet KRP personally, KRP's work remains a profound source of inspiration for him.
\end{quoting}

\begin{abstract}
    It is well known that free independence is equivalent to the vanishing of mixed free cumulants. The purpose of this short note is to build free products of $*$-probability spaces using this as the definition of freeness and relying on free cumulants instead of moments.
\end{abstract}

\section{Introduction}\label{sec:freeprodcons}
A non-commutative probability space (NCP) is a pair $(\A, \varphi)$ where $\A$ is an algebra over the complex numbers, with an unit say $\one_{\A}$, and $\varphi$ is a state, that is,  a linear functional on $\A$ such that $\varphi(\one_{\A})=1$. If $\A$ is a $*$-algebra and $\varphi$ is positive ($\varphi(a^*a)\geq 0$ for all $a\in \A$) then $(\A, \varphi)$ is called a $*$-probability space. In this case, it can be proved that $\varphi(a^*)=\overline{\varphi(a)}$ for all $a\in \A$ (see \mbox{\cite[p. 4]{spei}}).

We shall work with NCPs and $*$-probability spaces and often refer to them commonly as NCPs for convenience. Free independence is a central idea in non-commutative probability and is usually defined as follows.

\begin{defn}\label{def:freemoments} (Free independence)
Let $(\A, \varphi)$ be an NCP, and let $I$ be a fixed index set.
For each $i\in I$, let $\A_i$ be a unital sub-algebra of $\A$. Then
$(\A_i)_{i\in I}$ are called freely independent (or simply, free) if
$\varphi(a_1 \cdots a_k) = 0$
for every $k\geq 1$, whenever the following three conditions are satisfied:
\begin{enumerate}
    \item[(a)] for all $j$, \ $a_j\in\A_{i(j)}$ for some $i(j) \in I$,
    \item[(b)] $\varphi(a_j)=0$ for every $1\leq j\leq k$, and
    \item[(c)] neighbouring elements are from different sub-algebras: $i(1)\neq i(2)\neq \cdots \neq i(k-1)\neq i(k)$.
\end{enumerate}
We also say that variables $a_1, \ldots, a_n$ from $\A$ are freely independent (or simply, free) if the sub-algebras
generated by each of them are free.
\end{defn}

In Definition~\ref{def:freemoments}, if the $(\A_i)_{i\in I}$ are $*$-algebras, then they are usually said to be $*$-free (and likewise for variables $a_1, \ldots, a_n$ which are elements of a $*$-algebra). To keep things simple, we shall say ``free'' in either case.
It is easy to see that constants are free of any other variables.

Of course, one has to demonstrate the existence of freely independent variables before working with them. The classical analogue of this is the Kolmogorov construction which shows the existence of independent random variables.

Existence of freely independent variables $a_1, \ldots, a_n$ with given moments $\{\varphi(a_i^k), 1\leq i \leq n, k \geq 1\}$, or more generally, of freely independent NCPs $(\A_i, \varphi^{(i)})$, $i \in I$, is shown as follows (see, for example, \cite{spei}):
First a free product $\A$ of $(\A_i)_{i\in I}$ is constructed which contains copies of these algebras as sub-algebras. Then a state (moment functional) $\varphi$ is constructed on $\A$ such that for each $i \in I$, $\varphi$ agrees with $\varphi^{(i)}$ when restricted to $\A_i$, and the algebras $(\A_i)_{i \in I}$ become freely independent in $(\A, \varphi)$.

Just as classical independence translates to the vanishing of mixed cumulants, there is an equivalent, and arguably simpler, definition of freeness via the vanishing of mixed free cumulants (see Definition~\ref{def:alg}), which are free analogues of classical cumulants.

Our goal in this short note is to construct the free product of NCPs by relying solely on free cumulants and using the above vanishing property as the definition of freeness. We emphasise that our proof uses ideas from the free product construction outlined above. In particular, we first construct the free product of $(\A_i)_{i\in I}$ as outlined above. Then, instead of constructing the moment functional $\varphi$, we first develop free cumulants and arrive at the moment functional $\varphi$ on the free product space through these free cumulants.
\section{Free cumulants and free independence}
We need some background on free cumulants and moments.
Let $\NC(n)$, denote the set of all non-crossing partitions of $\{1, \ldots , n\}$  with the reverse refinement partial ordering. The symbols $\zero_n$ and $\one_n$ shall respectively denote the smallest and the largest partitions $\{\{1\}, \{2\}, \ldots ,\{n\}\}$ and $\{\{1, 2, \ldots , n\}\}$.  Note that $\NC(n)$ is a finite lattice. For any $\pi, \sigma\in \NC(n)$ let $\pi \vee \sigma$ denote the smallest element of $\NC(n)$ which is larger than both $\pi$ and $\sigma$. Being a finite lattice, it is equipped with a unique M\"{o}bius function, say $\mu_n$, on $\mathcal P_2(n)=\{(\pi, \sigma): \pi \leq \sigma, \pi, \sigma\in \NC(n)\}$. These functions extend naturally to a function, say $\mu$, on $\displaystyle{\sqcup_{n=1}^\infty \mathcal P_2(n)}$ such that the restriction of $\mu$ to $\mathcal{P}_2(n)$ is $\mu_n$. By an abuse of terminology, we shall refer to $\mu$ as the M\"{o}bius function on $\NC(n)$.

Let $(\mathcal{A},\varphi)$ be an NCP or a $*$-probability space. Suppose $\{a_i: i \in I \}$ is a collection of variables from $\mathcal{A}$. When $(\A, \varphi)$ is an NCP, the collection
\[
    \big\{\varphi\left(\Pi(a_i : i \in I)\right): \Pi\ \text{is a finite degree monomial}\big\}
\]
is called the collection of joint moments of $\{a_i : i \in I\}$. On the other hand, if $(\A, \varphi)$ is a $*$-probability space, then
\[
    \big\{\varphi\left(\Pi(a_i,a_i^{*}:i \in I)\right): \Pi\ \text{is a finite degree monomial}\big\}
\]
is the collection of joint $*$-moments of $\{a_i : i \in I\}$.

For every $n \geq 2$, we can extend $\varphi$ to a multilinear functional $\varphi_{n}$ on $\mathcal{A}^{n}$ as follows:
\begin{eqnarray} \label{eqn:mmt4.1}
\varphi_{n} (a_{1},a_{2},\ldots, a_{n}):= \varphi(a_{1}a_{2}\cdots a_{n}).
\end{eqnarray}
Then the multiplicative extension of $\{\varphi_n : n\geq 1\}$ to $\{\varphi_{\pi} : \pi \in \NC(n), n \geq 1\}$ is defined as follows:
If $\pi = \{V_{1}, V_{2}, \ldots, V_{r} \} \in \NC(n)$, then
\begin{eqnarray} \label{eqn:multi4.1}
\varphi_{\pi}[a_{1},a_{2},\ldots, a_{n}]:= \varphi(V_1)[a_{1},a_{2},\ldots,a_{n}]\cdots \varphi(V_r)[a_{1},a_{2},\ldots,a_{n}],
\end{eqnarray}
where
\[
    \varphi(V)[a_{1},a_{2},\ldots, a_{n}]:= \varphi_{s}(a_{i_1},a_{i_2},\ldots,a_{i_s})=\varphi(a_{i_1}a_{i_2}\cdots a_{i_s})
\]
for $V = \{i_1, \ldots, i_s\}$, where $i_1 < \cdots < i_s$. Note how the order of the variables has been preserved and the use of two types of braces $(\ )$ and $[\ ]$ in (\ref{eqn:mmt4.1}) and (\ref{eqn:multi4.1}). Also, in particular,
\begin{equation} \varphi_{\one_n}[a_1,a_2,\ldots, a_n] = \varphi_{n}(a_1,a_2,\ldots, a_n) = \varphi(a_1a_2\cdots a_n).
\end{equation}
We may now define free cumulants by using the M\"{o}bius function $\mu$.
\begin{defn} (Free cumulants) Suppose $(\A, \varphi)$ is an NCP. Then the \textit{joint free cumulant} functional $\kappa_n$ is defined on $\A^n$ as follows:
\begin{equation} \label{eqn:freecumdefn2}
\kappa_{n}(a_{1},a_{2},\ldots, a_{n}) := \sum_{\sigma \in \NC(n)} \varphi_{\sigma}[a_{1},a_{2},...,a_{n}]\mu(\sigma, \one_{n}), \ a_i\in \A, 1\leq i \leq n.
\end{equation}

Suppose $\{\A_i, i \in I\}$ are sub-algebras of $\A$. Then for any $n$, $\kappa_{n}(a_{1},a_{2},\ldots, a_{n})$, where $a_i\in \cup_{i\in I}\A_i$, is called a joint free cumulant. It is called a \textit{mixed free cumulant} if there are at least two variables $a_i, a_j$ such that $a_i\in \A_r, a_j \in \A_s, r, s \in I, \ r\neq s$.
\end{defn}
For any $n$, $\kappa_n (a^{\epsilon_1}, a^{\epsilon_2},\ldots, a^{\epsilon_n})$, where $\epsilon_i \in \{1, *\}, 1 \leq i \leq n$, is called a \textit{marginal free cumulant} of order $n$ of $\{a, a^{*}\}$. For a self-adjoint element $a$,
\[
    \kappa_{n}(a) := \kappa_n(a,a,\ldots, a)
\]
is called the $n$-th free cumulant of $a$.
Note that mixed/marginal free cumulants are all special cases of joint free cumulants. Clearly, the free cumulants $\kappa_n$  are also  multilinear functionals.

Just as in (\ref{eqn:multi4.1}), the family $\{\kappa_n: n \geq 1\}$ also has a multiplicative extension $\{\kappa_{\pi}: \pi \in \NC(n), n \ge 1\}$. Moreover, by making use of $\mu$ and its inverse, we can say that (see Proposition $11.4$ in \cite{spei}) for all $a_1, \ldots , a_n$,

\begin{align}
\kappa_{\pi}[a_1, a_2, \ldots, a_n] &= \sum_{\substack{\sigma \in \NC(n) \\ \sigma \leq \pi}} \varphi_{\sigma}[a_1, a_2, \ldots, a_n] \, \mu(\sigma, \pi)  \ \text{for all}\  \pi \in \NC(n), \ n \geq 1; \label{eqn: freecumdefn1} \\
\varphi_{\pi}[a_1, a_2, \ldots, a_n] &= \sum_{\substack{\sigma\in \NC(n) \\ \sigma \leq \pi}} \kappa_{\sigma}[a_1, a_2, \ldots, a_n] \ \ \text{for all}\ \ \pi \in \NC(n), \ n \geq 1; \label{eqn:momfreecum} \\
\varphi(a_1a_2\cdots a_n) &= \varphi_{\one_{n}}[a_1, a_2, \ldots, a_n] = \sum_{\substack{\sigma\in \NC(n) \\ \sigma \leq \one_n}} \kappa_{\sigma}[a_1, a_2, \ldots, a_n]. \label{eqn:momcum}
\end{align}

\noindent
Clearly, (\ref{eqn:freecumdefn2}) and (\ref{eqn:momcum}) provide a one-one correspondence between moments and free cumulants. In fact, similar formulas exist for moments and classical cumulants---the lattice of all partitions replaces $\NC(n)$ in these classical formulas. For an interesting discussion on the history of such classical moment-cumulant formulas, see the blog article \cite{speicherblog}. Anyway, since the moments and free cumulants are in one-one correspondence, it is natural that free independence may be defined via free cumulants.
\begin{defn}[Free independence]\label{def:alg}
Suppose $(\mathcal{A}, \varphi)$ is an NCP. Then sub-algebras $(\mathcal{A}_{i})_{i \in I}$ of $\mathcal{A}$ are said to be freely independent (or simply, free) if for all $n \geq 2$ and any $a_1, a_2,\ldots, a_n$ from $\cup_{i \in I} \mathcal{A}_i$, $\kappa_{n}(a_{1},a_{2},\ldots,a_{n}) = 0$ whenever at least two of the $a_i$'s  are from different $\mathcal{A}_j$'s. In particular, variables $a_1,
\ldots, a_n$ from $\A$ are said to be free if the sub-algebras generated by each of these variables are free.
\end{defn}
The following result is well known,
and first appeared in \cite{speimulti}.
See Theorem 11.16 of \cite{spei} for a proof. For a slightly shorter proof, see  Theorem 5.3.15 in \cite{anderson:guionnet:zeitouni}.

\begin{theorem}\label{thm:defequiv} Suppose $(\mathcal{A}, \varphi)$ is an NCP and $\mathcal{A}_{i}, i \in I$, are sub-algebras of $\mathcal{A}$. Then the following are equivalent. \vskip3pt

(a) The sub-algebras are free according to Definition~\ref{def:alg}.\vskip3pt

(b) The sub-algebras are free according to Definition~\ref{def:freemoments}.
\end{theorem}

\section{Free product of NCPs}
Before we state the main result, we recall the following standard construction (see, e.g., \cite[pp. 81-84]{spei}) of the \textit{free product algebra} $\mathcal{A}$ of the family of unital algebras $(\A_i)_{i \in I}$ corresponding to a family of NCPs $(\mathcal{A}_{i},\varphi^{(i)})_{i \in I}$. Using the universal properties for tensor products and direct sums, it may be shown that this construction satisfies the universal property for the free product (with identification of units) of unital algebras.
\vskip3pt

(i) First define the set of centered elements of $\mathcal{A}_{i}$ as
\[
    \mathring{\mathcal{A}}_{i} := \{a-\varphi^{(i)}(a)\one_{\mathcal{A}_{i}}:a\in \mathcal{A}_i\}, \ i \in I.
\]
The centered version of an element $a$ will be denoted by $\mathring{a}$.
\vskip3pt

(ii) Identify all the units $\one_{\mathcal{A}_{i}}$ as a single unit $\one_{\mathcal{A}} \equiv \one$.
\vskip3pt

(iii) Define the free product of $(\A_i)_{i \in I}$ as
\begin{equation}\label{eq:freeprodalg}
    \mathcal{A}=\mathbb{C}\one\ \displaystyle{\bigoplus_{k\geq 1}} \,\, \displaystyle{\bigoplus_{\substack{i_1\neq i_2 \neq \cdots \neq i_k \\ i_j \in I, \, 1 \leq j\leq k}}} \mathring{\mathcal{A}}_{i_{1}}\otimes \mathring{\mathcal{A}}_{i_{2}}\otimes\cdots \otimes \mathring{\mathcal{A}}_{i_{k}}.
\end{equation}
As an example, if $a_1\in \mathcal{A}_{1}$  and $a_2\in \mathcal{A}_{2}$, then the product $a_1a_2$ in $\mathcal{A}$ is defined as
\[
    \varphi^{(1)}(a_{1})\varphi^{(2)}(a_2) \, \one \, \oplus\varphi^{(1)}(a_1)\mathring{a}_2 \, \oplus \, \varphi^{(2)}(\mathring{a}_2)\mathring{a}_1 \, \oplus \, \mathring{a}_{1}\otimes \mathring{a}_{2}.
\]
We note here that each $\mathcal{A}_{i}$ has a `copy' in $\mathcal{A}$, namely $\mathbb{C} \one \oplus \mathring{\mathcal{A}}_i$. Also, the copies of $\mathring{\mathcal{A}}_{i}$, $i \in I$, inside $\A$ are pairwise disjoint.
\vskip3pt

(iv) With the observations in (iii), we can now define multiplication on $\mathcal{A}$ by concatenation and/or reduction as follows. Multiplication by scalars is defined in the obvious way. Suppose that $\Mixed{1}{k} \in \AlgMixed{1}{k}$ and $\Mixed[b][j]{1}{\ell} \in \AlgMixed[\mathcal{A}][j]{1}{\ell}$. Their product may now be defined inductively:
\begin{align*}
    (&\Mixed{1}{k}) (\Mixed[b][j]{1}{\ell}) \\
    &:= \begin{cases}
    \Mixed{1}{k} \otimes \Mixed[b][j]{1}{\ell}
    & \text{ if } i_k \ne j_1, \\
    \Mixed{1}{k - 1} \otimes \mathring{c} \otimes \Mixed[b][j]{2}{\ell} + \varphi_{i_k}(\mathring{c}) (\Mixed{1}{k - 1}) (\Mixed[b][j]{2}{\ell})  & \text{ if } i_k = j_1,
\end{cases}
\end{align*}
where $c = \mathring{a}_{i_k} \mathring{b}_{j_1}$ and the ``smaller order'' product $(\Mixed{1}{k - 1}) (\Mixed[b][j]{2}{\ell})$ is defined via the same rule. It is clear that $\mathcal{A}$ is an algebra.
\vskip3pt

(v) If the $\A_i$'s are $*$-algebras, then the $*$-operation on $\A$ is defined in the natural way via the specification $(a\otimes b)^*=b^*\otimes a^*$, and that makes $\mathcal{A}$ a $*$-algebra.

We are now ready to state and prove a precise result on the construction of the free product of NCPs based solely on free cumulants.
\begin{theorem}\label{theo:freeprodcons}
Let $(\mathcal{A}_{i},\varphi^{(i)})_{i \in I}$ be a family of NCPs (resp. $*$-probability spaces) and let $\mathcal{A}$ be as in \eqref{eq:freeprodalg}. Then there is a state $\varphi$ defined on $\mathcal{A}$ such that
\vskip3pt
(a) The sub-algebras (resp. $*$-sub-algebras) $(\mathcal{A}_{i})_{i \in I}$ are freely independent (resp. $*$-free independent) in $(\mathcal{A},\varphi)$.\vskip3pt

(b) The NCP $(\mathcal{A},\varphi)$ has the following universal property: Let $(\mathcal{B}, \psi)$ be an NCP such that for every $i \in I$, there exists a unital homomorphism $\Phi_i : \mathcal{A}_i \to \mathcal{B}$ such that  $\psi\circ\Phi_i=\varphi_i$,  and the images $(\Phi_i(\mathcal{A}_i))_{i\in I}$ are freely independent in $(\mathcal{B}, \psi)$. Then there exists a homomorphism $\Phi$ between $(\mathcal{A}, \varphi)$ and $(\mathcal{B}, \psi)$, uniquely determined, such that $\Phi |_{\mathcal{A}_i} = \Phi_i$ for every $i\in I$.
\end{theorem}

The rest of this note is devoted to proving Theorem~\ref{theo:freeprodcons}. We will construct the state $\varphi$ by first constructing the free cumulant functionals $\kappa_n$, $n \ge 1$, in such a way that all mixed free cumulants vanish. This will prove Part (a). Part (b) is the same as Part (2) of Proposition 6.6 in \cite{spei} and we refer the reader to \cite{spei} for its proof.

We break the proof down into three steps.
\vskip3pt

\noindent
\textbf{Step 1. Construction of the free cumulant functionals $\boldsymbol{\kappa_n}$.}
Let $\kappa_n^{(l)}$ denote the free cumulant functionals corresponding to  $(\A_l, \varphi^{(l)}), \ l\in I$. We shall now define free cumulant functionals $\kappa_n: \A^n \rightarrow \mathbb{C}$, $n \geq 1$, making sure that $\kappa_n$ agrees with $\kappa_n^{(l)}$ on each $\A_l^{n}$, $ l \in I$ (as per our earlier discussion, $\mathcal{A}_l^n$ is identified with its copy inside $\mathcal{A}^n$ for every $l$).
Moreover, all mixed free cumulants must vanish. Thus, we start by defining, for every $n$,
\begin{equation} \label{eq:freecumstep1}
    \kappa_n(a_1, \ldots, a_n):= \begin{cases} \kappa_n^{(l)}(a_1, \ldots, a_n) &\text{if all } a_i\text{'s are from the same sub-algebra } \mathcal{A}_{l} \text{ for some } l, \\
0 &\text{if at least two } a_i\text{'s are from different sub-algebras}.
\end{cases}
\end{equation}

By construction, these functionals $\kappa_n$ are multilinear on $\A_{l_1}\times\cdots  \times \A_{l_n}$ for every choice of $l_1, \ldots, l_n$.

Call $\{a_i\}, \ 1\leq i \leq n$, ``pure'', if each of them belong to one of the sub-algebras $\mathcal{A}_l$.
For pure $\{a_i\}$, define $\kappa_{\pi}$ by the multiplicative extension of $\kappa_n$ defined above as
\begin{equation}\label{eq:freepure}
    \kappa_{\pi}[a_1, \ldots, a_n] := \prod_{V=\{i_1, \ldots,  i_s\}\in \pi}\kappa_{|V|}(a_{i_{1}}, \ldots, a_{i_{s}}),
\end{equation}
where the right side is defined by (\ref{eq:freecumstep1}).

We shall now use (\ref{eq:freecumstep1}) and (\ref{eq:freepure})
to define all free cumulants. \vskip3pt

First consider monomials in pure elements, say $b_1, \ldots, b_m$. For simplicity, write
\begin{eqnarray*}
b_1&=&a_1\cdots a_{s_{1}},\\
b_j &=& a_{s_{j - 1} + 1} \cdots a_{{s_{j}}}, \ j \geq 2.
\end{eqnarray*}
That is, $b_j$ is the product of $s_j - s_{j - 1}$ pure elements. Note that here we always work with the obvious \textit{minimal} representation of a product of pure variables in order to write $b_j$'s in a unique way (up to scalar multiples), i.e. we stipulate that consecutive variables in the descriptions of the $b_j$'s are from different sub-algebras. In view of the relation that free cumulants of elements from free products must obey (see Equation (11.11) in Theorem 11.12 of \cite{spei}), to define $\kappa_m(b_1, \ldots , b_m)$, we have no other choice but the following:
\begin{equation}\label{eq:freemixed}
    \kappa_m(b_1, \ldots, b_m) := \sum_{\substack{\pi \in \NC(s_m) \\ \pi \vee \sigma = \mathbf{1}_{s_m}}} \kappa_\pi[a_{1}, \ldots, a_{s_{1}}, \ldots, a_{s_{m - 1} + 1}, \ldots, a_{s_{m}} ],
\end{equation}
where $\sigma = \{ \{1, \ldots, s_1\}, \ldots, \{s_{m - 1} + 1, \ldots, s_m \}\} \in \NC(s_m)$. Note that each $\kappa_\pi$ on the right hand side of \eqref{eq:freemixed} has already been defined by \eqref{eq:freepure}. Further, because of the minimal representation used in the descriptions of the $b_j$'s, it is clear that $\kappa_m$ is well-defined. Now for arbitrary $\pi \in \NC(m)$ we define
\begin{equation}\label{eq:freemixedfull}
    \kappa_\pi[b_1, \ldots, b_m] :=\prod_{V=\{i_1, \ldots,  i_s\}\in \pi}\kappa_{|V|}(b_{i_{1}}, \ldots, b_{i_{s}}).
\end{equation}
Note that each term on the right hand side of \eqref{eq:freemixedfull} has already been defined in \eqref{eq:freemixed}. \begin{remark}
We note here that the earliest combinatorial description of the relations between moments and free cumulants was given in \cite{speimulti}, which was later fully developed in \cite{krwawcykspei}. A special case of formula (\ref{eq:freemixed}) already appeared in \cite{speimulti}.
\end{remark}

Now, by multilinearity, we extend this definition to $\kappa_{\pi}[b_1, \ldots, b_m]$ where the $b_j$'s are polynomials instead of monomials. In particular, the functionals $\kappa_n:\mathcal{A}^{n}\to \mathbb{C}$ are now completely defined for all $n$ and they agree with $\kappa_n^{(l)}$ on $\A_l^n$ for all $l \in I$.
If the above construction is carefully followed, it is also clear that $\kappa_1(a^*) = \overline{\kappa_1(a)}$ for all $a\in \A$, when $\A$ is a $*$-algebra. \vskip3pt

\noindent
\textbf{Step 2. Construction of the state $\boldsymbol{\varphi}$ on $\boldsymbol{\A}$.} Having defined $\{\kappa_n : n \ge 1\}$, we define, as in (\ref{eqn:momfreecum}),
\begin{equation}\label{eq:momentfreecum}\varphi_{\pi}[a_1, a_2, \ldots, a_n] := \sum_{\substack{\sigma\in \NC(n)\\ \sigma \leq \pi}} \kappa_{\sigma}[a_1,a_2,\ldots, a_n] \ \ \text{for all}\ \ \pi \in \NC(n), \ n \geq 1.
\end{equation}
In particular,
\[
    \varphi(a):=\varphi_\one[a]=\kappa_\one[a]=\kappa_1(a).
\]
Clearly, $\varphi$ is a state since it is linear and $\varphi(\one)=\kappa_1(\one)=1$. Hence  $(\mathcal{A}, \varphi)$ is an NCP.
Moreover, by construction,
\[
    \varphi(a)=\kappa_1(a)=\kappa_1^{(l)}(a)=\varphi^{(l)}(a)\ \ \text{for all}\ \ a\in \mathcal{A}_{l}, \ l \in I.
\]
Note that the $\kappa_n$'s are indeed the free cumulant functionals corresponding to $\varphi$ due to the relation (\ref{eq:momentfreecum}).

Now since all mixed free cumulants vanish (see (\ref{eq:freecumstep1})),  the algebras $\mathcal{A}_i$, as sub-algebras of
$\mathcal{A}$, are free with respect to the state $\varphi$. This proves Part (a) of Theorem~\ref{theo:freeprodcons} for NCPs which are not $*$-probability spaces. Incidentally, due to the manner of our construction of the free cumulants, and the one-one correspondence between moments and free cumulants given in (\ref{eqn:freecumdefn2}) and (\ref{eqn:momcum}), it is immediate that $\varphi$ is uniquely defined.
\vskip3pt

\noindent
\textbf{Step 3. Positivity of $\boldsymbol{\varphi}$ for $\boldsymbol{*}$-algebras.}
For $*$-algebras, we have to additionally show that  $\varphi(a^*a) \ge 0$ for all $a \in \mathcal{A}$. By the relation \eqref{eqn:momcum} between moments and free cumulants, we have, for all $a \in \mathcal{A}$, that
\[
    \varphi(a^*a) = \kappa_2(a^*, a) + \kappa_1(a^*)\kappa_1(a).
\]
But, the second term on the right side above is non-negative since $\kappa_1(a^*) = \overline{\kappa_1(a)}$ as noted before.
Hence we only need to check that the ``variance'' $\kappa_2(a^*, a) \ge 0$.
We now prove three lemmas which will accomplish this.

\begin{lemma}\label{lem:constfree}
    For any $a = \Mixed{1}{k} \in \AlgMixed{1}{k}$, we have $\kappa_2(\mathbf{1}, a) = \kappa_2(a, \mathbf{1}) = 0$.
\end{lemma}
\begin{proof}
By \eqref{eq:freemixed}, we have
\begin{align*}
    \kappa_2(\mathbf{1}, a) &= \kappa_2(\mathbf{1}, \Mixed{1}{k}) \\
                         &= \sum_{\substack{\pi \in \NC(k + 1) \\ \pi \vee \sigma = \mathbf{1}_{k + 1}}} \kappa_\pi [\mathbf{1}, \Pure{i_1}, \ldots, \Pure{i_k}],
\end{align*}
where $\sigma = \{\{1\}, \{2, \ldots, k + 1\}\}$. Note that each $\pi$ in the above sum must have $1$ in a block $V_{\pi, 1}$ of length $\ge 2$. Suppose $V_{\pi, 1} = \{1, j_1, \ldots, j_t\}$. If $i_{j_1} = \cdots = i_{j_t}$, then
\[
    \kappa_{|V_{\pi, 1}|}(\mathbf{1}, \Pure{i_{j_1}}, \ldots, \Pure{i_{j_t}}) = \kappa_{|V_{\pi, 1}|}^{(i_{j_1})}(\mathbf{1}, \Pure{i_{j_1}}, \ldots, \Pure{i_{j_t}}) = 0,
\]
because constants are free of non-constant elements in the original sub-algebras. Otherwise, $t > 1$ and there exist $1 \le t' <  t'' \le t$ such that $i_{j_{t'}} \ne i_{j_{t''}}$. Then by \eqref{eq:freecumstep1}, $\kappa_{|V_{\pi, 1}|}(1, \Pure{i_{j_1}}, \ldots, \Pure{i_{j_t}})$ is zero. It follows that $\kappa_{\pi}[1, \Pure{i_1}, \ldots, \Pure{i_{s_i}}] = 0$ for each $\pi$ such that $\pi \vee \sigma = \mathbf{1}_{k + 1}$. Thus $\kappa_2(\mathbf{1}, a) = 0$, and, by a similar argument, so is $\kappa_2(a, \mathbf{1})$.
\end{proof}

\begin{lemma}\label{lem:variance}
Let $a = \Mixed{1}{k} \in \AlgMixed{1}{k}$ and $b = \Mixed[b][j]{1}{\ell} \in \AlgMixed[\mathcal{A}][j]{1}{\ell}$. Then $\kappa_2(a^*, b)$ is non-zero only if $k=\ell$ and $i_u = j_u$ for all $1 \le u \le k$, in which case
\begin{equation}\label{eq:freemixedtwo}
    \kappa_2(a^*, b) = \prod_{u = 1}^k \kappa_2(\PureStar{i_u},\Pure[b]{j_u}).
\end{equation}
In particular, $\kappa_2(a^*, a) \ge 0$.
\end{lemma}
\begin{proof}

We have, by \eqref{eq:freemixed}, that
\begin{equation}\label{eq:freemixedtwo-firststep}
    \kappa_2(a^*, b) = \sum_{\stackrel{\pi \in \NC(k + \ell)}{\pi \vee \sigma = \mathbf{1}_{k + \ell}}} \kappa_\pi [\PureStar{i_k}, \ldots, \PureStar{i_1}, \Pure[b]{j_1}, \ldots, \Pure[b]{j_{\ell}}],
\end{equation}
where $\sigma = \{ \{1, \ldots, k\}, \{k + 1, \ldots, k + \ell\} \}.$ Suppose that $\kappa_2(a^*, b) \ne 0$. Consider a partition $\pi$ appearing in \eqref{eq:freemixedtwo-firststep} for which the corresponding $\kappa_\pi \ne 0$. Clearly, $\pi$ cannot be $\mathbf{1}_{k + \ell}$, for otherwise $\kappa_\pi$ would be zero by centeredness. Take two indices $i, j \in \{1, \ldots, k+ \ell\}$ that are in the same block, with the smallest possible value for $j - i > 0$.

Now, by freeness, we must have that the variables corresponding to $i, j$ come from the same algebra in order for $\kappa_\pi$ to be non-zero. If $j - i > 1$, then, there cannot exist another index $i' \in \{i + 1, \ldots, j - 1\}$ which belongs to the same partition block as $i$ and $j$, for then it would violate the minimality of $j - i$. For the same reason, no two $i', j' \in \{i + 1, \ldots, j - 1\}$ can belong to the same block of $\pi$. We conclude that all entries in $\{i + 1, \ldots, j - 1\}$ must constitute singleton blocks. But this would force $\kappa_{\pi}$ to be zero by centeredness. It follows that we must have $j = i + 1$. Now, since consecutive elements of both $a$ and $b$ come from different algebras, in order for $\kappa_{\pi}$ to be non-zero, $i, j$ cannot both correspond to elements of $a$, or of $b$---the only possibility is that $i = k$ and $j = k + 1$ and that $i_1 = j_1$.

Now, there cannot exist an $\overline{i} < k$ (resp. $\overline{j} > k+1$) such that $\overline{i}, k$ (resp. $k+1, \overline{j}$) are in the same block of $\pi$, for otherwise there would exist an $\overline{i}$ (resp. $\overline{j}$) with $k - \overline{i}$ (resp. $\overline{j} - (k+1)$) being the smallest, and then by the argument of the above paragraph, $k = \overline{i} + 1$ (resp. $\overline{j} = k+1 + 1$). This would make $\kappa_{\pi}$ zero as $\overline{i}, k$ both correspond to elements of $a$ (resp. $k+1, \overline{j}$ to elements of $b$). This means that $\{k, k+1\} $ is a block of $\pi$, and hence
\[
    \kappa_{\pi} [\PureStar{i_k}, \ldots, \PureStar{i_1}, \Pure[b]{j_1}, \ldots, \Pure[b]{j_\ell}] = \kappa_2(\PureStar{i_1}, \Pure[b]{j_1}) \,\, \kappa_{\pi \setminus \{k, k + 1\} }[\PureStar{i_k}, \ldots, \PureStar{i_2}, \Pure[b]{j_2}, \ldots, \Pure[b]{j_\ell}].
\]
Noting that $\pi \setminus \{k, k + 1\} \in \NC(k + \ell - 2)$, we can apply the same argument to conclude that $i_2 = j_2$ and $\{k - 1, k + 2\}$ is a partition block of $\pi$ and so on. It is therefore clear that we must have $k = \ell$ and $i_u = j_u$ for all $1 \le u \le k$ and that only the partition $\pi_{2k} = \{ \{1, 2k\}, \{2, 2k - 1\}, \ldots, \{k, k + 1\} \}$ has a non-zero $\kappa_\pi$. Since $\pi_{2k}$ does appear in the sum in \eqref{eq:freemixedtwo-firststep}, \eqref{eq:freemixedtwo} follows. The last assertion follows because
\[
    \kappa_2(a^*, a) = \prod_{u = 1}^k \kappa_2(\PureStar{i_u}, \Pure{i_u}) = \prod_{u = 1}^k \kappa_2^{(i_u)}(\PureStar{i_u}, \Pure{i_u}),
\]
and the variances $\kappa_2^{(i_u)}(\PureStar{i_u}, \Pure{i_u})$ are all non-negative.
\end{proof}
Now we handle elements $a$ which are sums of elements of a particular $\AlgMixed{1}{k}$.
\begin{lemma}\label{lem:samealg}
    Let $a = \sum_{t = 1}^r a_t$, where $a_t = \Pure{t, i_1} \otimes \cdots \otimes \Pure{t, i_k} \in \AlgMixed{1}{k}$. Then $\kappa_2(a^*, a) \ge 0$.
\end{lemma}
\begin{proof}
    By multilinearity and Lemma~\ref{lem:variance},
\begin{align*}
    \kappa_2(a^*, a) &= \sum_{1 \le s, t \le r} \kappa_2(a_s^*, a_t) \\
                     &= \sum_{1 \le s, t \le r} \prod_{u = 1}^k \kappa_2(\PureStar{s, i_u}, \Pure{t, i_u}) \\
                     &= \mathbf{1}^\top M \mathbf{1},
\end{align*}
where $M$ is an $r \times r$ matrix with entries $M_{st} = \prod_{u = 1}^k \kappa_2(\PureStar{s, i_u}, \Pure{t, i_u})$. It suffices, therefore, to show that $M$ is positive semi-definite. To that end, note that $M$ can be represented as
\[
    M = M_{1} \odot \cdots \odot M_{k},
\]
where $(M_u)_{st} = \kappa_2(\PureStar{s, i_u}, \Pure{t, i_u})$ and $\odot$ denotes elementwise product of matrices (also known as the Schur-Hadamard product). The Schur product theorem \cite{schur1911bemerkungen} states that the elementwise product of positive semi-definite matrices is positive semi-definite. Therefore we will be done if we show that $M_u$ is positive semi-definite. This follows easily because, for any $x\in \C^r$,
\begin{align*}
    x^* M_u x &= \sum_{s, t} \bar{x}_s x_t \kappa_2 (\PureStar{s, i_u}, \Pure{t, i_u}) \\
              &= \kappa_2\Big(\big(\sum_s x_s \Pure{s, i_u}\big)^*, \sum_s x_s \Pure{s, i_u} \Big) \\
              &= \kappa_2(w^*, w),
\end{align*}
where $w = \sum_s x_s \Pure{s, i_u} \in \Alg{i_u}$, and $\kappa_2$ agrees with $\kappa_2^{(i_u)}$ on $\Alg{i_u}$.
\end{proof}
Now we consider a general element $a \in \mathcal{A}$, which can be expressed, for some $r \ge 0$, as
\[
    a = c\mathbf{1} \oplus \bigoplus_{i = 1}^r a_{i},
\]
where $c \in \mathbb{C}$, each $a_i$ is a sum of some elements in $\AlgMixed{1}{s_i}$, i.e.
\[
    a_{i} = \sum_{j = 1}^{t_i} a_{ij}, \ a_{ij} \in \AlgMixed{1}{s_i}, 1\le j \le t_i, 1 \le i \le r,
\]
and the hosts $\AlgMixed{1}{s_i}$ are all different. We expand $\kappa_2$ by multilinearity as
\begin{equation} \label{eq:gen-decomp}
    \kappa_2(a^*, a) = |c|^2 \kappa_2(\mathbf{1}, \mathbf{1}) + \bar{c}\sum_{i} \kappa_2(\mathbf{1}, a_{i}) + c \sum_{i}\kappa_2(a_{i}^*, \mathbf{1}) + \sum_{i} \kappa_2(a_i^*, a_i) +\sum_{i\ne i'} \kappa_2(a_{i}^*, a_{i'}).
\end{equation}
Note that $\kappa_2(\mathbf{1}, \mathbf{1}) = 0$, as $\kappa_2$ agrees with $\kappa_2^{(l)}$ for $l \in  I$. By Lemmas~\ref{lem:constfree} and \ref{lem:variance} we see that all but the third sum in \eqref{eq:gen-decomp} are zero, and the third sum itself is non-negative by Lemma~\ref{lem:samealg}.

This completes the proof of positivity of $\varphi$ for $*$-algebras. The proof of Theorem~\ref{theo:freeprodcons} is now complete.

\section*{Acknowledgements}
We thank Roland Speicher for pointing us to crucial references and to the history of the classical moment-cumulant formulas. The detailed comments from the anonymous referee have led to a substantial improvement in the article.

\bibliographystyle{alpha}
\bibliography{bib.bib}

\end{document}